\documentclass[[letterpaper, 10 pt, conference]{ieeeconf}
\IEEEoverridecommandlockouts
\usepackage{cite}
\usepackage{amsmath,amssymb,amsfonts}
\usepackage{graphicx}
\usepackage{textcomp}
\usepackage{verbatim}
\usepackage{caption}
\usepackage{subcaption}
\usepackage{multirow}
\def\BibTeX{{\rm B\kern-.05em{\sc i\kern-.025em b}\kern-.08em
    T\kern-.1667em\lower.7ex\hbox{E}\kern-.125emX}}

\title{\LARGE \bf Electric Vehicle -
Power Grid Incorporation Using Distributed Resource
Allocation Approach}

\author{J. Moyalan, M. Sawant, Bhagyashree U., A. Sheikh, S. Wagh and N. Singh\\EED, VJTI, Mumbai, India
\thanks{J. Moyalan, M. Sawant, Bhagyashree U., A. Sheikh, S. Wagh and N. Singh are with the Electrical Engineering Department, Veermata Jijabai Technological Institute, Mumbai 400019, INDIA
        {\tt\small jrmoyalan\_m17@ee.vjti.ac.in}}%
}

\begin{document}

\maketitle

\begin{abstract}
Models of multi-agent systems can be found all around us. One of the objectives of multi-agent systems is to define local control laws in order to achieve a desired global state of the system. This paper utilises the concept of Distributed Resource Allocation (DRA) approach for successful incorporation of a large number of Plug-in Electric Vehicles (PEVs) with the power grid. DRA approach is implemented using the concept of achieving output consensus. A fixed number of PEVs are considered which are connected to the grid for a certain time interval. The PEV batteries can be charged as well discharged during this time interval. Charging and discharging of PEVs from the grid is further divided into time slots and are coined as strategies. The goal of this paper is to obtain a well-inclined charging strategy for each  PEV connected to the grid such that it satisfies the power grid objective in terms of the smoothening factor of grid load profile. 
\end{abstract}

\section{Introduction}
Generation of power has always come with big price to be paid by the environment. However, in the present era where everything is run by electricity, it is difficult to lessen this burden on nature on a large scale. As a result, energy management is one of the biggest concerns of the present electrical grid systems. However, this problem can be culminated on a small scale through various methods. Of these, the most innovative method is the incorporation of plug-in Electric Vehicle (PEV) with the power grid, \cite{wu2011review,shao2011demand,pollet2012current}. The presence of several PEVs in the grid offer both rewards and challenges in functioning of an efficient electricity grid \cite{6195074, 6933935}. PEVs act like virtual batteries which can give as well as store electricity which helps in dealing with the sudden changes in the load demand which is difficult to be met by the generators in short time without significant cost. On the other hand, if several PEVs are connected to a single power grid then the optimal utilisation of each PEVs keeping in mind the constraints of the PEVs itself, network and that of the generators is cumbersome to achieve. 
Several approaches have been employed to realise the imminent task of optimal involvement of PEVs with the power grid. They can be broadly classified into centralised and decentralised approaches. In centralised approach, like in \cite{clement2010impact,sundstrom2012flexible}, a central body collects information from all the PEVs and solve an optimisation problem and based on the solution sends strategies to each participating agent of the system for desired output. However, there are several drawbacks to these methods such as the inability to account the owners choices \cite{gan2013optimal,richardson2012local}. In decentralised approach, there is no such thing as a central body. The whole system is divided into subparts and each part solves an optimisation problem based on the information obtained from the remaining subparts of the system as given in \cite{he2012optimal,rezaei2014packetized}. One such method is the concept of Distributed Resource Allocation (DRA) using output concensus where agents or subparts of the system make decisions based on the local information and coordinate with each other within the system to achieve a desirable global state.
\begin{table}[h]
\caption{Symbol Description}
\begin{tabular} {|c|l|}
\hline
\textbf {Symbol} & \textbf { Description} \\
\hline
 $\mathcal{S}$ & Set of nodes of a multi-agent system\\
\hline
 $\mathcal{L}$  & Set of edges connecting the nodes of a multi-agent system\\
\hline
\multirow{3}{*}{$\mathcal{A}$} & A nonnegative matrix whose elements satisfy the\\
    &  following: $a_{kj}=1$ if $(k,j) \in  \mathcal{L}$; $a_{kj}=0$ if $(k,j) \notin  \mathcal{L}$ \\
\hline
$\mathcal{N}_k$ & Set of neighbours of node $k$\\
\hline
$A_k$ & Total active power supply of grid at $k^{th}$ time slot\\
\hline
 $R_k$ & Total reactive power supply of grid at $k^{th}$ time slot\\
\hline
$\beta$ & Barrier function\\
\hline
\multirow{3}{*}{$x^i_k$} & Active power charging strategy at $k^{th}$ time step\\
 & of $i^{th}$ PEV\\
\hline
\multirow{3}{*}{$y^i_k$} & Reactive power charging strategy at $k^{th}$ time step\\
 & of $i^{th}$ PEV\\
\hline
$K^i$ & Number of time steps allotted to $i^{th}$ PEV\\
\hline
$soc^i_K$ & Desired state of charge of $i^{th}$ PEV\\ 
\hline
$soc^i_o$ & Initial state of charge of $i^{th}$ PEV\\
\hline
$\overline{soc}^i_o$ & Upper limit of $i^{th}$ PEV charger\\
\hline
$\underline{soc}^i_o$ & Lower limit of $i^{th}$ PEV charger\\
\hline
$t^i_k$ & Time width of $k^{th}$ time step of $i^{th}$ PEV\\
\hline
$\overline{p}^i$ & Nominal power of $i^{th}$ PEV charger\\
\hline
$\mu$ & Commitment factor controlled by the PEV owner\\
\hline
$\eta$ & Smoothing factor controlled by the power grid manager\\
\hline
$s^i$ & Auxillary slack variable of $i^{th}$ PEV\\
\hline
$Q^i$ & Total reactive power available of $i^{th}$ PEV\\
\hline
$\overline{q}^i_k$ & Available reactive power of of $i^{th}$ PEV at $k^{th}$ time step\\
\hline
\end{tabular}
\end{table}
The remainder of this paper is partitioned as follows: Section II introduces the problem statement and its analogies with suggested approach. Section III presents the prerequisite to understand the working of the proposed approach for load management problem. Section IV gives a detailed description of DRA approach using output consensus. Section V presents the proposed approach for multiple PEVs incorporation with power grid based on the explanations of previous sections. Section VI presents case studies showing suggested approach under the influence of single and multiple PEVs. In Section VII some conclusions are given.

\section{Problem Statement and Analogies}

A distribution system transformer is considered which supply power to two types of customers: commercial and transient. Commercial customers include industrial regions and occupational areas whose energy load profile do not change with time. The latter type of customers have their load profile varying with time. In this paper, PEVs are considered as transient customers which use the power from the grid to charge their batteries upto the maximum level.
\begin{figure}[htbp]
\centerline{\includegraphics[width=\linewidth]{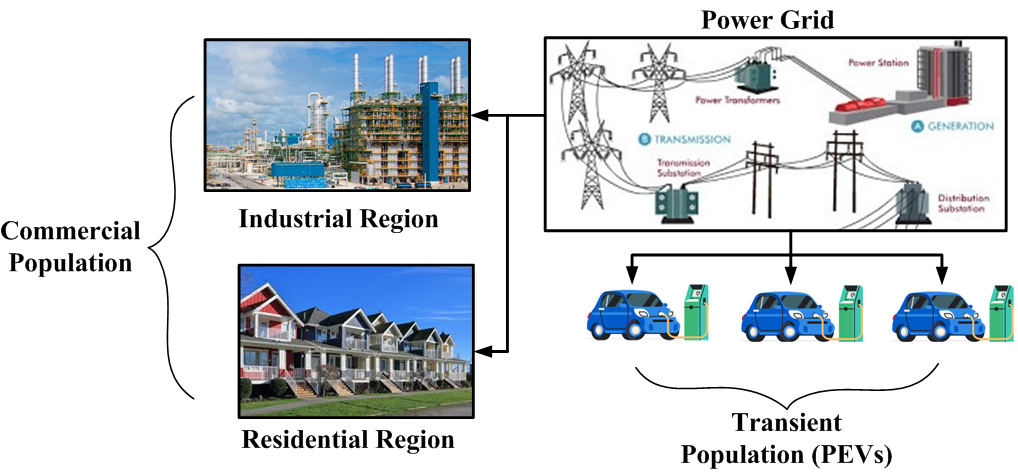}}
\caption{Illustrative diagram of power grid connected to commercial population and transient population.}
\label{fig}
\end{figure}
For smooth incorporation of PEVs with the grid, different payoff functions are proposed at different time instants depending on the electricity demand. Accordingly, there are $k$ time slots for corresponding $k$ different payoff functions. The active and reactive power given by the grid to its customer is represented by array $A$ and $R$, respectively.
\begin{equation}
A_{k}=a_{k}+\sum_{i=1}^{N}x^{i}_{k}/t^i_k,  \quad               
R_{k}=r_{k}+\sum_{i=1}^{N}y^{i}_{k}               
\end{equation}
where $a_{k}$ and $r_{k}$ represents active and reactive power respectively delivered by the distribution transformer at $k^{th}$ time slot in the absence of any PEVs. Similarly $x^i_{k}$ and $y^i_{k}$ represents active energy and reactive power of $i^{th}$ PEV at time slot $k$. Each PEV represents two types of population namely real power (W) and reactive power (VAr).

The analogies proposed for energy populations can be illustrated in terms of strategies \cite{ovalle2017escort}. Depending on the values of $k$, there are $k$ strategies where every strategy provides a certain payoff for individuals settled on it. The PEVs have the option of discharging their battery at certain time slots in order to satisfy its objectives. The discharged energy will be supplied to the commercial population by the grid and provides portion of the payoff of that time slot to that particular PEV. In this way, transient population (PEVs) can force commercial population to mitigate to other strategies. However, the total number of populations spanning across all the strategies remains constant.

Since different time slots have different payoff functions imposed by the grid, the aim of this paper is to find the energy consumed by each PEV at different time slots. The final energy profile of all PEVs should be such that it results in several beneficiary services to the grid such as smoothening of grid output profile and minimizing of load shifting. 

\section{Preliminaries}
\subsection{Graph Theory}\label{AA}

In the multi-agent system considered in this paper, modeling of the communication network with the help of graph allows the agents to coordinate their decisions as given in \cite{godsil2001g}. A graph is described by the triplet $\mathcal{C}=( \mathcal{S},\mathcal{L},\mathcal{A})$. $\mathcal{S}={\{1,. . ., K}\}$ represents the set of nodes, $\mathcal{L}\subseteq\mathcal{S}\times\mathcal{S}$ represents the set of edges connecting the nodes and matrix $\mathcal{A}$ represents a $K\times K$ nonnegative matrix whose elements satisfies: $a_{kj} = 1$ if and only if $(k,j)\in\mathcal{L}$, and $a_{kj}=0$ if and only if $(k,j)\notin\mathcal{L}$. The nodes and edges of the graph corresponds to agents and communication channels of the multi-agent system respectively. The neighbours of node $k$ i.e. the set of nodes that are able to receive/send information from/to node $k$ is represented by $\mathcal{N}_{k}={\{j \in \mathcal{S} : (k,j) \in \mathcal{L}}\}$.

Graphical modeling of the multi-agent system can be completed under the following assumptions for simplicity:

1) $a_{kk}=0 \quad \forall k \in \mathcal{S}$ i.e. no nodes are connected to itself.

2) $a_{kj}=a_{jk}$ i.e. channels are bidierctional.

The graph laplacian matrix of $\mathcal{C}$ can be defined as $K \times K$ matrix $L(\mathcal{C}) = [l_{kj}]$ as follows:

\begin{equation} 
l_{kj} = 
\begin{cases}
\sum_{j \in \mathcal{S}} a_{kj}, \quad if \quad k = j\\
-a_{kj}, \quad \quad if \quad  k \neq j\\
\end{cases}
\end{equation}
\subsection{System Passivity Framework}

Consider a dynamical system represented by (\ref{eq:equation1})
\begin{equation}
\begin{aligned}
\dot{x} = f(x,u)\\
y = h(x,u)\\
\end{aligned}
\label{eq:equation1}
\end{equation}

where $f:R^n \times R^p \rightarrow R^n$ is locally Lipschitz, $h: R^n \times R^p \rightarrow R^p$ is continuous, $f(0,0) = 0$, and $h(0,0) = 0$. The system has the same number of inputs and outputs.

The system given in (3) is said to be passive if there exists a continuously differentiable positive semidefinite function $V(x)$ (called storage function) such that
\begin{equation}
u^Ty \ge \dot V = \frac{\partial V}{\partial x}{f(x,u)}, \quad \forall (x,u) \in R^n \times R^p
\end{equation}

Moreover, it is said to be

\begin{itemize}
\item lossless if $u^Ty= \dot V \quad \forall (x,u) \in R^n \times R^p$.
\item strictly passive if $u^Ty \ge \dot V + \psi (x)$ for some positive definite function $\psi$, and for $\forall (x,u) \in R^n \times R^p$.
\end{itemize}

\subsection{Barrier Function Formulation}

Barrier function $\beta (x)$ as described in  \cite{bravo2015distributed} is used as a power constraint to prevent the dynamical equations of resource allocation problem from violating its constraints. Consider a system with $(a,b)$ as constraint set where $x$ is the input variable. Then $\beta (x) $ has the  following properties:

\begin{itemize}
\item $\beta (x) $ is monotonically increasing continuous function defined in $(a,b)$.
\item $\beta (x) \rightarrow -\infty$, when $x \rightarrow a$.
\item $\beta (x) \rightarrow \infty$, when $x \rightarrow b$. 
\end{itemize}

Barrier function can be thought as the derivative of a convex function which prevents the control signal from going outside its feasible domain.

\section{Features of DRA}

Consider a multi-agent system of $n$ agents connected by a communication network  that is characterized by the weighted graph $\mathcal{C}=( \mathcal{S},\mathcal{L},\mathcal{A})$. The dynamical model of the system can be represented by the following differential equations:

\begin{equation}
\Gamma^S_k :
\begin{cases}
\dot x_k = f(x) \\

y_k = g(x)\\
\end{cases}
\end{equation}


where $\Gamma^S_k$ represents the system as a whole, $y_k \in R$ denotes the output of subsystem $k$ and $x_k$ is the state of subsystem $k$.

The objective of all agents is to drive the system to a desired global state which corresponds to the desired grid objectives. Consider the situation where each agent has information of its output and that of its neighbours', i.e. the $k^{th}$ agent knows the value of $y_k$ and that of $y_j$ for all  $j \in \mathcal{N}_k$. The available information is used by the agents to formulate their control law which drives their $\dot x$. This is shown in (6):
\begin{equation}
\Gamma^C_k : \dot x_k = u_k(y_k,y_j), \forall j \in \mathcal{N}_k
\label{dot_eq}
\end{equation}
\subsection{Control Objective}

In multi-agent optimization setting, the desired global state mostly corresponds to a common consensus value. Such problems are called output consensus problems. The definition of output consensus as defined in \cite{bravo2015distributed} is given below.

\textit{Definition 1:} Consider the set of subsystems given in (5) and (6). It can be said that output consensus is reached if $lim_{t \rightarrow \infty} \mid y_k(t) - y_j(t) \mid = 0$, for all $i, j = 1, . . ., n$ where $y_k(t)$ is the output of the subsystem $i$ at time t.

The case of reaching output consensus is achieved under the constraint applied to the state variables as follows:
\begin{equation}
\sum_{k=1}^{n} x_k^i = X,
\end{equation}
where $X \in R$ is the sum of all the state values in the given time frame. In the present problem statement $x_k^i$ is the electrical energy transferred between the PEVs and the grid and $X$ represents the total amount of electrical energy required to completely charge the battery of each PEV.

The control objective of the multi-agent system can be summarized as follows:

\textbf{1.} Satisfying the constraint (7).

\textbf{2.} Driving equation (5) to output consensus.

\subsection{Resource Allocation Dynamics}\label{SCM}

In order to achieve the desired global state, local control laws $u_1, u_2, . . ., u_n$ needs to be designed to be applied to the multi-agent system which is proposed as follows: 
\begin{equation}
u_k(y_k,y_j) = \sum_{j \in \mathcal{N}_k} a_{kj} (y_j - y_k); \forall k = 1,...,K, \forall j \in \mathcal{N}_k
\end{equation}
It is easily seen that (8) satisfies our constraint (7) if following conditions are met:

\begin{itemize}
\item $ \sum^{K}_{k=1} x_k^i(0) = X$\\
\item $\sum^{K}_{k=1} \dot x_k^i = 0$
\end{itemize}

\subsection{Convergence to Output Consensus}
The feedback interconnection outlook is appropriate for the multi-agent system represented by (5) utilizing the control law (6) as displayed in Fig. 2. The equilibrium point $(x^*)$ of the feedback interconnection after application of proposed DRA dynamical equation (8) must satisfy the following statement which is adapted from \cite{bravo2015distributed}:

\begin{figure}[htbp]
\centerline{\includegraphics[scale=0.5]{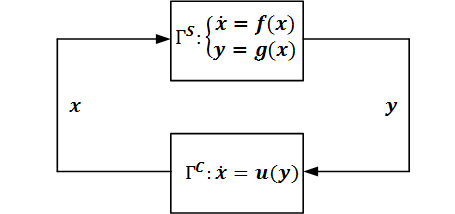}}
\caption{Feedback interconnection of systems (5) and (6).}
\label{fig}
\end{figure}

\textit{Statement 1:} Consider the feedback interconnection shown in Fig. 2. having its equilibrium point at $x^*$ and let the steady state output of $\Gamma^S$ be given by  $y^*=g(x^*)$. If $u(y_k,y_j) \quad \forall j \in \mathcal{N}_k$ is given by (8) and the communication graph $\mathcal{C}$ is connected, then $y^*_k=y^*_j \quad \forall k,j = 1,..., K$ where $y^*_k$ is the $k^{th}$ element of the vector $y^*$.

Using the definition of output consensus problem, Statement 1 states that the output consensus will be obtained if the equilibrium point $x^*$ is asymptotically stable. To check the stability of $x^*$ the dynamics of $\Gamma^S$ and $\Gamma^C$ is expressed in error coordinates.

The system $\Gamma^S$ is written in error coordinates as follows,
\begin{equation}
\Gamma^S_e:
\begin{cases}
\dot e_x = f^e(e_x)\\
e_y = g^e(e_x)\\ 
\end{cases}
\end{equation}
where $e_x = x-x^*$ and $e_y=y-y^*$. Also $f^e(e_x)=f(x)$ and $g^e(e_x)=g(x)-g(x^*)$ for all $x \in R^n$. Since $x^*$ is an equilibrium point of (5), it can be seen that $f^e(0)=0$, and $g^e(0)=0$ .

Following assumption is made on (9):

\textit{Assumption 1:} Consider the dynamical system (9). If $f^e(0)=0$, then $e_x=0$.

The \textit{Assumption 1} guarantees the existence of unique rest points for (5).

The dynamics of (6) that implements (8) is also expressed in error coordinates using Laplacian of $ \mathcal{C}$ as follows:
\begin{equation}
\dot e_x=-L ( \mathcal{C}) e_y
\end{equation}
Now consider the following Statement given in \cite{bravo2015distributed} reformulated as follows:

\textit{Statement 2:} The multi-agent system expressed in error coordinates given by (8) is passive and lossless from the input $e_y$ to the output $-e_x$, if  $x(0)$ and $x^*$ satisfies the resource constraint (7), i.e., $\sum^K_{k=1}x^*_k=X$ and $ \mathcal{C}$ is connected.

The concept of passivity can be explored along with Statement 2 to validate the stability of equilibrium points of (5),\cite{bravo2015distributed}. The property that the feedback interconnection of two passive systems generally results in stable rest points is utilized to ensure output consensus is achieved under the configuration as shown in Fig. 2. The following theorem adapted from  \cite{bravo2015distributed} is used to summarize the requirements to reach output consensus.

\textit{Theorem 2:} Consider the feedback interconnection of system (5) and (6) having its equilibrium point at $x^*$ where (8) defines the $u(y)$. Following conditions are assumed:

\begin{enumerate}
\item[{i)}]The graph $\mathcal{C}$ of the system (6) is connected.

\item[{ii)}] The resource constraint (7) is satisfied by $x^*$ and $x(0)$.

\item[{iii)}] Assumption 1 is satisfied by the system (5) expressed in error coordinates with respect to $x^*$. Moreover it is strictly passive from the input $e_x$ to the output $e_y$ with radially unbounded storage function.
\end{enumerate}

Then (5) reaches output consensus. 

\section{DRA Application: PEV Load Management}

The application of incorporation of PEVs with a microgrid is described in this paper based on (8) and the properties of output consensus problem.

\subsection{Constraints on PEV variables}

The energy constraints considered are similar to that given in \cite{ovalle2017escort} in the application of PEV integration with the grid.
\begin{equation}
\sum^{K^i}_{k=1}  x^i_{k} = soc^i_K - soc^i_0,
\end{equation}
\begin{equation}
x^i_{k} \le \overline{soc}^i - (soc^i_0 + \sum^{\Omega}_{\omega=1}  x^i_{\omega} - x^i_{k}), 
\end{equation}
\begin{equation}
x^i_{k} \ge \underline{soc}^i - (soc^i_0 + \sum^{\Omega}_{\omega=1}  x^i_{\omega} - x^i_{k}),
\end{equation}
$$\forall \Omega = {\{1,2,..., K^i}\}, \forall k = {\{1,2,...,\Omega}\},$$
\begin{equation}
-t^i_k \overline{p}^i \le x^i_{k} \le t^i_k \overline{p}^i, \quad \forall k = {\{1,2,...,K^i}\}
\end{equation}
where $soc^i_K$ represents as desired state of charge (SoC) (in Watt-hour) at the end of time window, $soc^i_0$ represents initial SoC (in Watt-hour), $t^i_k$ is the length of $k^{th}$ time step (in hours), and $\overline{p}^i$ is the nominal power of the charger. Constraint (11) is equivalent to the state variable constraint given in (7). Constraints (12) and (13) defines the accumulated SoC at a particular time instant which is not allowed to cross upper limit $ \overline{soc}^i$ as well as lower limit $ \underline{soc}^i$. Constraint (14) defines limits of energy consumption and injection rates by PEV depending upon the limits of charger as well as duration of time steps.

The following reactive power constraints are used:
\begin{equation}
\sum^{K^i}_{k=1} y^i_{k} + s^i = 0,
\end{equation}
\begin{equation}
- \overline{q}^i_{k} \le y^i_{k} \le \overline{q}^i_{k}, \quad \forall k = {\{1,2,...,K^i}\},
\end{equation}
\begin{equation}
-Q^i \le s^i \le Q^i.
\end{equation}
where $y^i_{k}$ defines total reactive power of charger at each time step, $s^i$ defines the auxillary slack variables and available reactive power $ \overline{q}^i_{k}$ is defined as
\begin{equation}
\overline{q}^i_{k} = \pm \sqrt{( \overline{p}^i)^2 - (x^i_{k}/t^i_k)^2}
\end{equation}
Also, the total reactive power $Q^i$ available in the time window is defined as
\begin{equation}
Q^i = \sum^{K^i}_{k=1} \sqrt{( \overline{p}^i)^2 - (x^i_{k}/t^i_k)^2}
\end{equation}
Equations (15) and (17) are defined in such a way that sum of all the absolute values of contributions of reactive power in the time window has to be less or equal to $Q^i$.
\begin{equation}
-Q^i \le \sum_{k=1}^{K^i}y_k^i \le Q^i
\end{equation}
and the portion of $Q^i$ that is not used in the time window are assigned to the slack variables.

\subsection{Formulation of Output Functions}
Payoff function of each time step would be the output function of this particular resource allocation problem. The inclusion of commitment factor $\mu^i$ and a smoothing factor $\eta $ correlates objectives of the power grid as well as that of PEV owners. $\mu^i$ is controlled by the PEV owners according to their level of satisfaction while, $\eta$, a parameter controlled by the power grid manager monitors the sudden changes in the active and reactive power. Therefore,
\begin{equation}
\underline{\mu} \le \mu^i < 1, \quad 0 \le \eta \le 1
\end{equation}
where $\underline \mu$ is the minimum allowed limit of commitment.

The payoff function for active power strategies are
\begin{equation}
f^i_{k}(x^i_{k}) = -(1-\mu)(x^i_{k}-x^{i*}_{k})/t^i_k - \mu \eta A_{k}
\end{equation}
$$\quad \quad \quad - \mu(1-\eta)(2A_{k}-A_{k-1}-A_{k+1}),$$
where $x^{i*}_{k}$ is the PEV owner's preferred charging rate at the $k^{th}$ time step  and $\mu$ is the mean value of all $\mu^i$.

Similarly, the payoff function for Reactive power strategies are defined as follows:
\begin{equation}
g^i_{k}(y^i_{k}) = -(1-\mu)y^i_{k} - \mu \eta R_{k}
\end{equation}
$$\quad \quad \quad \quad \quad \quad - \mu(1-\eta)(2R_{k}-R_{k-1}-R_{k+1}),$$

Depending upon the values of $\mu$ and $\eta$ it is possible to modify the output profile of total as well as single phase active and reactive power profiles. When $\mu = 0$, flattening/smoothening objectives are ignored, while payoff functions give importance to local references of load distribution. When $\eta = 0$ and $\mu > 0$, smoothing objective is prioritize while flattening objectives are neglected. Finally, when $\eta = 1$ and $\mu > 0$ flattening objectives are preferred to smoothing objective. However, in each above cases, the utility grid manager has direct or indirect control over these parameters which is used as a trade-off between grid objectives as wll as economic and social benefits to PEV owners.

The inclusion of the constraints given in (11)-(14) for active power strategies and in (15)-(17) and (20) for reactive power strategies in the output function is implemented by a simple modification of adding the barrier function proposed in Section III-C to the payoff function given in (22) and (23). Therefore the modified output functions for both active and reactive power strategies is given as follows:
\begin{equation}
l^i_k(x^i_k) = f^i_{k}(x^i_{k}) + \epsilon \beta (x^i_k)
\end{equation}
\begin{equation}
p^i_k(y^i_k) = g^i_{k}(y^i_{k}) + \epsilon \beta (y^i_k)
\end{equation}
where $\epsilon>0$ is a small positive constant. The value of $\epsilon$ is decided in such a way that the effect of $\beta$ is truncated when the control input is away from the boundary values. For active power strategy, $x^i_k \in (a,b)$ where $a \in max[ \underline{soc}^i - (soc^i_0 + \sum^{\Omega}_{\omega=1}  x^i_{\omega} - x^i_{k}),-t^i_k \overline{p}^i]$ and $b \in min[ \overline{soc}^i - (soc^i_0 + \sum^{\Omega}_{\omega=1}  x^i_{\omega} - x^i_{k}), t^i_k \overline{p}^i]$. Similarly, for reactive power strategy, $y^i_k \in (a,b)$ where $a = - \overline{q^i_k}$ and $b = \overline{q^i_k}$.

The functionality of barrier formulation (11)-(14) is to constrain the active power strategies in (15)-(17) and reactive power in (20). For example, consider the active power strategy $x^i_k$ to be close to the upper bound. This corresponds to the higher value of $\beta (x^i_k)$ as it is monotonically increasing function of $x^i_k$, which ultimately results in higher payoff function value (24). $\dot x^i_k$ becomes negative (according to (6) and (8)) when the above condition occurs and therefore the value of $x^i_k$ will decrease without violating the upper.

\subsection{Connectivity of Proposed Problem}

The connectivity of the strategies as mentioned in Section III-A is necessary for output consensus to achieve. The connectivity graph of the strategies proposed in this paper is governed by (1), (22) and (23). An example of the connectivity graph is shown in Fig. 3. for three PEVs connected to the grid, i.e., $i={\{1,2,3}\}$ and grid providing payoff functions for four time slots i.e., $k={\{1,2,3,4}\}$ where $(k,i)$ represents the strategy with $k^{th}$ time slot for $i^{th}$ PEV.
\begin{figure}[h!]
\centerline{\includegraphics[width=\linewidth]{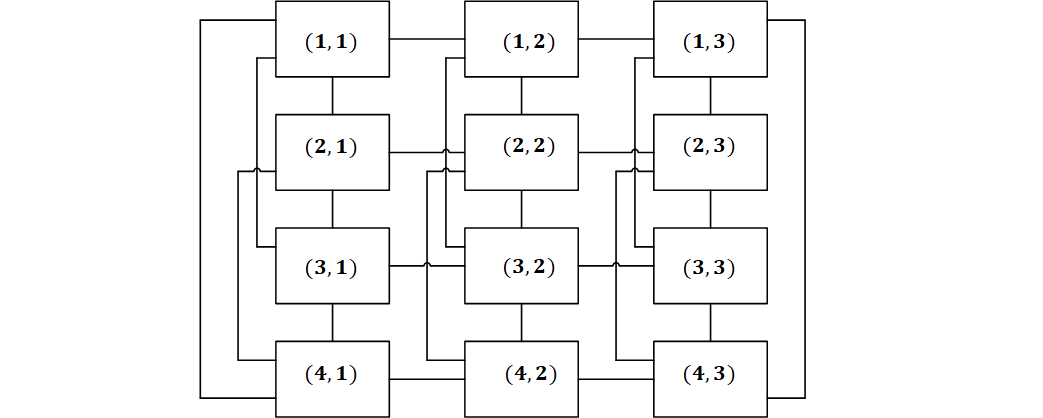}}
\caption{System partial information connectivity diagram }
\label{fig}
\end{figure}
\section{Results}

For the implementation scenario, all PEVs are considered with $soc_0=15.5 kWh$ and $soc_K=16 kWh$. Therefore, each PEV requires $500Wh$ to reach its $soc_K$. Consider $\overline{soc^i}$ and $\underline{soc^i}$ for all PEV batteries to be $18.5 kWh$ and $13.5 kWh$ respectively. Let the nominal power of PEV charger to be $3kW$. Let there be $K=10$ strategies or time slots for each PEV and each time slot has a width of 1 hour.

\begin{figure*}[ht]

\begin{subfigure}[b]{0.5\linewidth}
\centering
\includegraphics[width=\linewidth]{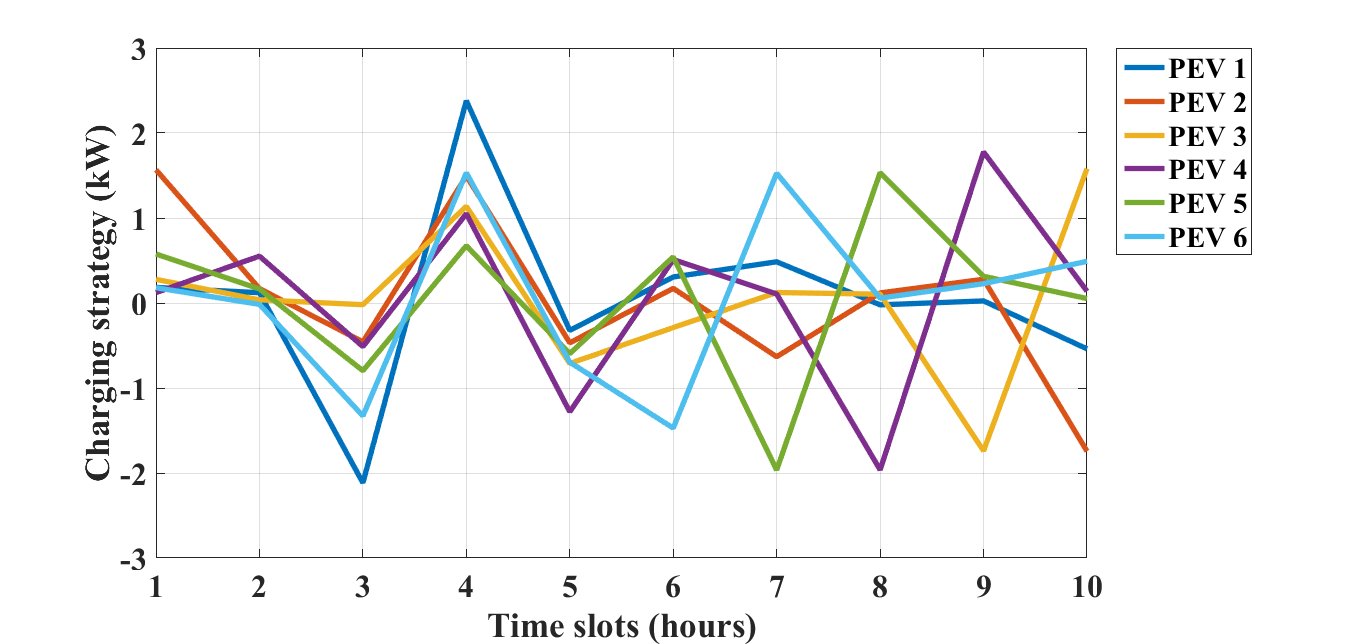}
\caption{Charging strategies for six PEVs in given time slots}
\end{subfigure}
\begin{subfigure}[b]{0.5\linewidth}
\centering
\includegraphics[width=\linewidth]{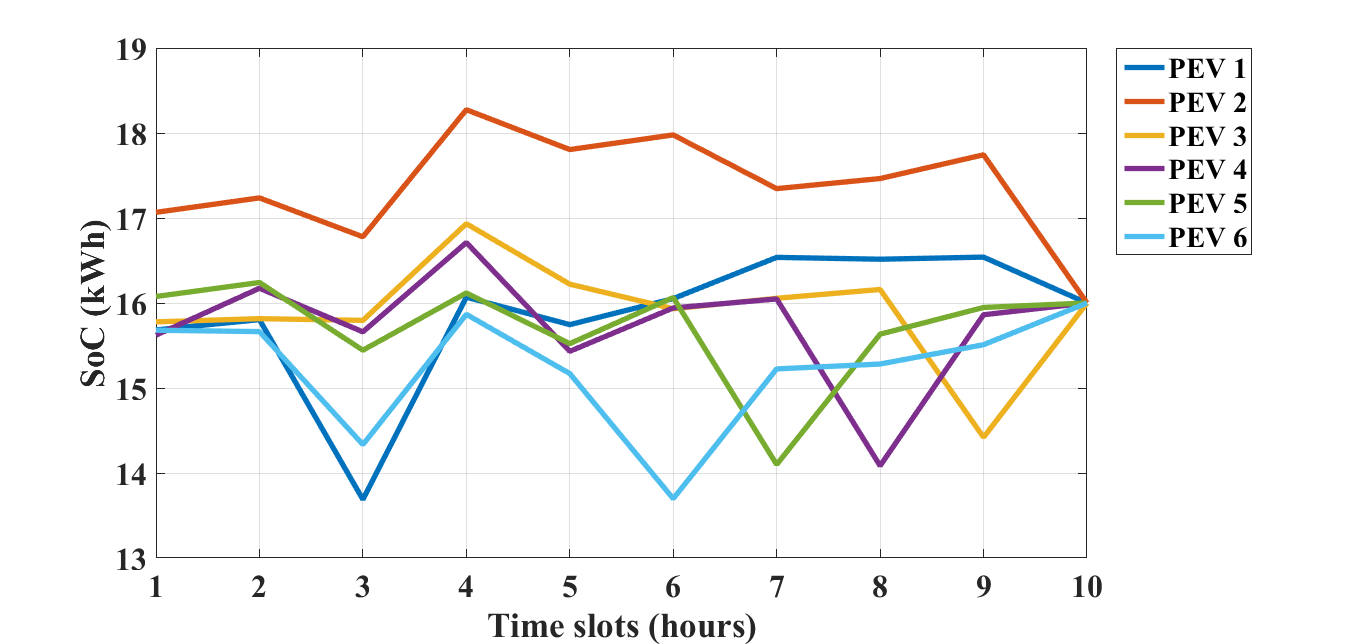}
\caption{State of charge of batteries of six PEVs in given time slots}
\end{subfigure}
\begin{subfigure}[b]{0.5\linewidth}
\centering
\includegraphics[width=\linewidth]{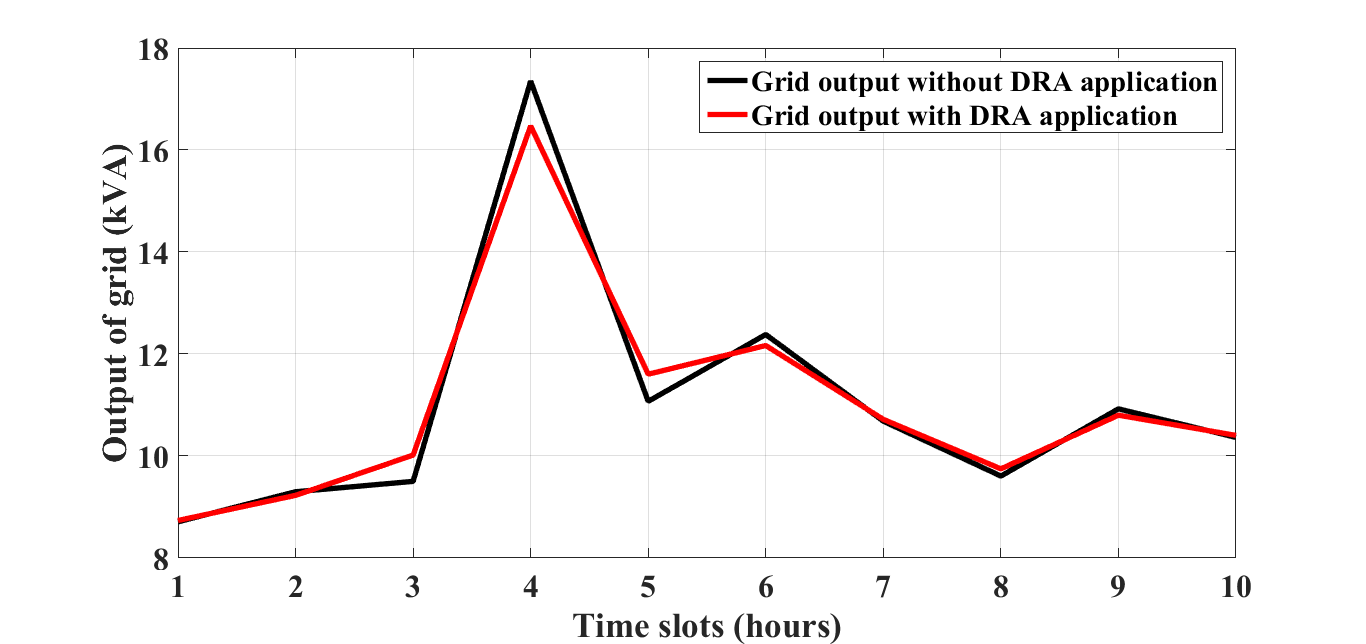}
\caption{Output of grid connected with single PEV}
\end{subfigure}
\begin{subfigure}[b]{0.5\linewidth}
\centering
\includegraphics[width=\linewidth]{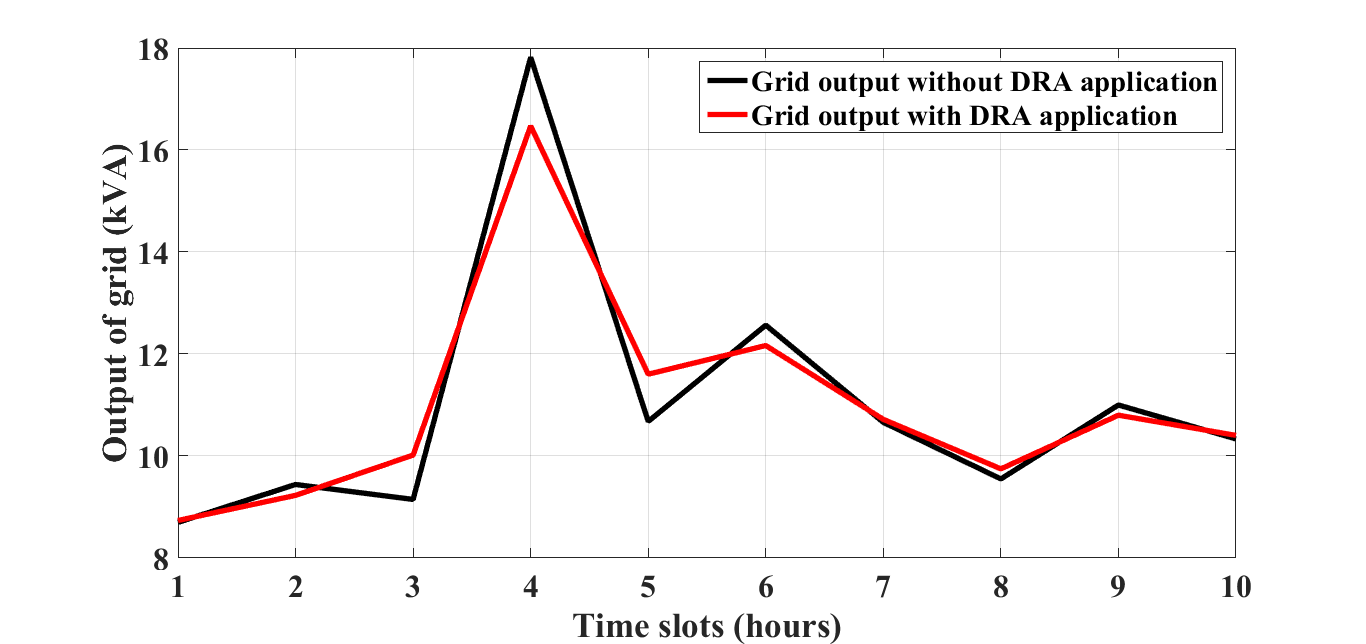}
\caption{Output of grid connected with six PEVs}
\end{subfigure}
\caption{Six PEVs connected to a power grid sharing similar conditions.}
\end{figure*}

\subsection{Effect of number of PEVs}

Let the external parameters be fixed at values $\mu = 0.5$ and $\eta = 0.8$. Fig. 4 displays the results of incorporation of six PEVs with the power grid using DRA approach of attaining output consensus of payoff function. Fig. 4a displays the charging strategy of the six PEVs at each time slot, in which negative value corresponds to the discharging instance of respective PEVs. The output consensus reached using DRA approach helps in achieving the conditions described in the output function. The payoff function (22) and (23) penalize the deviations between desired and proposed charging rate while taking into account the smoothing of grid output load profile. As a result, at consensus, PEV discharging its battery at certain time slots corresponds to the tradeoff between the grid objectives and economic interest of vehicle owners. Such occasion can arise when there is sudden increment in load profile of commercial customers, which then is tried to be minimized utilizing the energy present in the PEV batteries. Fig. 4b displays the SoC of PEV batteries after each time slot which converges to $16 kWh$ thereby achieving the objective of the vehicle owners in the given time slots. Fig. 4c shows the output of the power grid in the given time frame with and without the implementation of DRA approach using a single PEV. It can be seen that the output values after implementing DRA method is slightly more smooth compared to the output graph of the grid plotted without using DRA method. Moreover, as the number of PEVs participating in energy transfer process increases, the smoothening effect of output load profile becomes more evident as seen in Fig. 4d. Therefore, the following inference can be reached: more the number of PEVs, more apparent power available for the smoothening of the output load profile of the grid.

\subsection{Effect of $\mu$ and $\eta$}
\begin{figure}[htbp]
\centerline{\includegraphics[width=\linewidth]{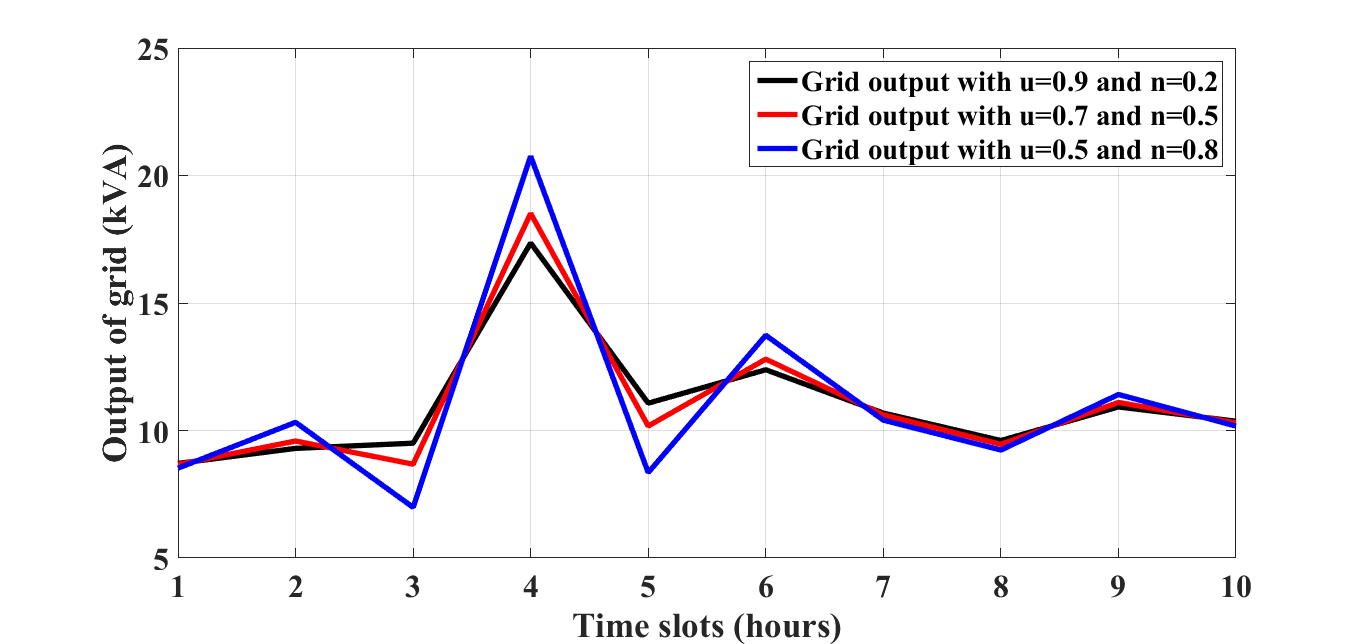}}
\caption{Power grid output for different values of $\mu$ and $\eta$}
\label{fig}
\end{figure}
As described in Section V-B, the commitment of the PEV to the grid for the given time frame is parametrized by $\mu$ while the smoothing factor $\eta$ characterizes the degree to which the output load profile of the grid is smoothened. It can be interpreted from Fig. 5 that for $\mu>0$, the smoothing objective of the power grid output depends on the value of $\eta$ in an inverse fashion. As the value of $\eta $ increases in the interval $[0,1]$, the importance of smoothing objective goes on decreasing. When $\mu=0$, then the payoff function (22) and (23) gives value to the  local references of load distribution. Therefore, the value of $\mu$ defines the closeness of actual charging profile of each PEV with that of their desired charging profile in the given time frame.

\section{Conclusion}
Distributed load management of power grid is done utilising the concept of output consensus, an application of DRA. Variation in the output load of the grid with respect to commitment factor $\mu$ and smoothening factor $\eta$ is also analysed. A case study has been conducted whose results concluded that the proposed approach gives desirable performances in PEV load management with respect to smoothening of grid supply. The future aim of this paper is to conduct similar real time case studies where the number of PEVs connected to a power grid varies for a given time frame. 

\bibliography{References}
\bibliographystyle{IEEEtran}
\end{document}